\nonstopmode
\documentclass[11pt]{jdg-p}
\usepackage{amssymb,amscd}
\usepackage{enumerate}
\usepackage[numbers]{natbib}
\usepackage[all]{xy}
\newdir{ >}{{}*!/-10pt/@{>}}
\newdir^{ (}{{}*!/-5pt/@^{(}}
\newdir_{ (}{{}*!/-5pt/@_{(}}
\def\cgaps#1{}
\def\Cgaps#1{}
\allowdisplaybreaks

\def\undersetbrace#1\to#2{\underbrace{#2}_{#1}}								 
\def\oversetbrace#1\to#2{\overbrace{#2}^{#1}}
\def\AMSunderset#1\to#2{\underset{#1}{#2}}
\def\AMSoverset#1\to#2{\overset{#1}{#2}}

\swapnumbers

\newtheorem*{prop*}{Proposition}

\newtheorem*{thm*}{Theorem}

\newtheorem*{lem*}{Lemma}

\newtheorem*{rem*}{Remark}
\newtheorem*{cor*}{Corollary}

\numberwithin{equation}{subsection}

\def\ign#1{}             

\def\X{\mathfrak X}
\def\al{\alpha}
\def\be{\beta}
\def\ga{\gamma}
\def\de{\delta}

\def\si{\sigma}
\def\ta{\tau}

\def\Ga{\Gamma}
\def\De{\Delta}

\def\Ph{\Phi}

\def\Om{\Omega}
\def\i{^{-1}}
\def\x{\times}
\def\p{\partial}
\let\on=\operatorname
\def\L{\mathcal L}

\def\AMSonly#1{}

\def\R{\mathbb R}

\def\Tr{\on{Tr}}
\def\vol{\on{vol}}
\def\Vol{\on{Vol}}

\def\g{\bar{g}}

\def\Ric{\on{Ricci}}
\def\Scal{\on{Scal}}

\def\vp{\vphantom{\big(\big)}}
\def\Met{{\mathcal M}}
\def\na{\nabla}
\def\wt{\widetilde}
\def\N{\mathcal N}
\def\NN{\R_{>0}\,g_0}
\begin{document}

\title[Metrics on the  manifold of all Riemannian metrics]
{Sobolev metrics on the  manifold of all Riemannian metrics}
\author{Martin Bauer, Philipp Harms, Peter W. Michor}
\address{
Martin Bauer, Peter W. Michor:
Fakult\"at f\"ur Mathematik, Universit\"at Wien,
Nordbergstrasse 15, A-1090 Wien, Austria.
\newline\indent
Philipp Harms: Edlabs, Harvard University, 44 Brattle street, Cambridge, MA 02138
}
\email{Bauer.Martin@univie.ac.at}
\email{pharms@edlabs.harvard.edu}
\email{Peter.Michor@univie.ac.at}
\thanks{All authors were supported by 
FWF Project 21030, MB was supported by FWF Project P24625}
\date{\today}
\keywords{Space of all Riemannian metrics, Sobolev metric}
\subjclass{58D17, 58E30, 35A01}

\begin{abstract}
On the manifold $\Met(M)$ of all Riemannian metrics on a compact manifold $M$ one can 
consider the natural $L^2$-metric as described first by \cite{Ebin70}.
In this paper we consider variants of this metric which in general are of higher order.
We derive the geodesic equations, we show that they are well-posed under some conditions and induce 
a locally diffeomorphic geodesic exponential mapping. We give a condition when Ricci flow is a 
gradient flow for one of these metrics. 
\end{abstract}

\maketitle

\section{Introduction}\label{nmb:1}

On the manifold $\Met(M)$ of all Riemannian metrics on a compact manifold $M$ one can 
consider the natural $L^2$-metric. 
It was first described by \cite{Ebin70}. Geodesics and curvature on it were described 
by \cite{FreedGroisser89} and \cite{Gil-MedranoMichor91} who also described the Jacobi fields and 
the exponential mapping. This was extended to the space of non-degenerate bilinear structures on $M$
in \cite{Gil-MedranoMichorNeuwirther92} and restricted to the space of almost Hermitian structures 
in \cite{Gil-MedranoMichor94}. In his thesis \cite{Clarke09a} which was published in two subsequent papers
\cite{Clarke09b,Clarke09c}, Brian Clarke showed that geodesic 
distance for the $L^2$-metric is a positive topological metric on $\Met(M)$, and he determined 
the metric completion of $\Met(M)$. In contrast, 
it was shown in \cite{Michor98, Michor102} that the natural $L^2$-metric on the space 
of immersions from a compact manifold into a Riemannian manifold has 
indeed vanishing geodesic distance. This also holds for the right invariant $L^2$-metric on 
diffeomorphism groups \cite{Michor102}, and even on the Virasoro-Bott group \cite{Michor122} where
the geodesic equation is the KdV-equation.

In this paper, guided by the results of \cite{Michor118, Michor119, Michor120}, we investigate 
stronger metrics on $\Met(M)$ than the $L^2$-metric. These are metrics of the following 
form:
\begin{align*}
G_g(h,k)&=\Phi(\Vol)\int_M g^0_2(h,k) \vol(g) &&\text{see \ref{nmb:4.2}}
\\
\text{or  }&=\int_M \Phi(\Scal).g^0_2(h,k) \vol(g)  &&\text{see \ref{nmb:4.3}}
\\
\text{or  }&=\int_M g^0_2((1+\De)^ph,k) \vol(g)  &&\text{see \ref{nmb:4.4}}
\end{align*}
where $\Ph$ is a suitable real-valued function, $\Vol=\int_M \vol(g)$ is the total volume of 
$(M,g)$, $\Scal$ is the scalar curvature of $(M,g)$, and where $g^0_2$ is the induced metric on 
$\binom{0}{2}$-tensors.
We describe all these metrics uniformly as 
$$
G^P_g(h,k)=\int_M g^0_2(P_gh,k) \vol(g)=\int_M \Tr(g\i. P_g(h). g\i. k) \vol(g),
$$ 
where 
$P_g:\Ga(S^2T^*M) \to \Ga(S^2T^*M)$
is a positive, symmetric, bijective pseudo-differential operator of order $2p, p\geq 0,$ 
depending smoothly on the metric $g$. 
We derive the geodesic equation for the general metric and all particular cases. 
We show that under certain assumptions on $P_g$ the geodesic equation is well posed and that the 
geodesic exponential mapping is a diffeomorphism from a neighborhood of the 0 section in the 
tangent bundle $T\Met(M)$ onto a neighborhood of the diagonal in $\Met(M)\x \Met(M)$.
The assumptions are satisfied for the metrics in \ref{nmb:4.2} and \ref{nmb:4.4}, but not for the metric 
in \ref{nmb:4.3}. In many cases the curve $(1-t)g_0$ can be reparameterized as a geodesic. In each 
case we can estimate its length, getting conclusions about geodesic incompleteness.

Finally we derive a condition on $P_g$ which is sufficient for the Ricci vector field to be a 
gradient field in the $G^P$-metric.

We thank the referee for very helpful remarks.

\section{Notation}\label{nmb:2}
\subsection{Metric on tensor spaces}\label{nmb:2.1}

A Riemannian metric $g:TM \x_M TM \to \R$ will equivalently be interpreted as 
$$\flat = g: TM \to T^*M \quad \text{and} \quad \sharp=g\i: T^*M \to TM.$$ 
The metric $g$ can be extended to the cotangent bundle $T^*M=T^0_1M$ by setting
$$g\i(\alpha,\beta)=g^0_1(\alpha,\beta)=\alpha(\beta^\sharp)$$
for $\alpha,\beta \in T^*M$, 
and the product metric 
$$g^r_s = \bigotimes^r g \otimes \bigotimes^s g\i$$
extends $g$ to all tensor spaces $T^r_s M$. 
A useful formula is
\begin{equation*}
g^0_2(h,k)= \on{Tr}(g\i h g\i k) \quad \text{for $h,k \in T^0_2M$ if $h$ or $k$ is symmetric.}
\end{equation*}
For a proof using orthonormal frames see \cite{Michor119}.
In this work, traces always contract the first two free appropriate tensor slots:
$$\on{Tr}: T^r_s M \to T^{r-1}_{s-1}, \quad \on{Tr}^g: T^r_s M \to T^{r-2}_s M.$$

\subsection{Directional derivatives of functions}\label{nmb:2.2}

We will use the following ways to denote directional derivatives of functions, in particular in 
infinite dimensions.
Given a function $F(x,y)$ for instance,
we will write:
$$ D_{(x,h)}F \text{ or } dF(x)(h) \text{ as shorthand for } \partial_t|_0 F(x+th,y).$$
Here $(x,h)$ in the subscript denotes the tangent vector with foot point $x$ and direction $h$. 
If $F$ takes values in some linear space, we will identify this linear space and its tangent space. 
We use calculus in infinite dimensions as explained in \cite{MichorG}.

\subsection{Volume density}\label{nmb:2.3}

The \emph{volume density} on $M$ induced by the metric $g$ is given by
$\on{vol}(g)=\on{vol}(g)\in\Gamma(\vol(M)),$ where $\vol(M)$ denotes the volume bundle.
The \emph{volume} of the manifold with respect to the metric $g$ is given by
$\on{Vol}=\int_M \on{vol}(g).$
The integral is well-defined since $M$ is compact. If $M$ is oriented we may identify the volume 
density with a differential form.  
Furthermore we have the following formula for the first variation of the volume density
(see for example \cite[Section 3.6]{Michor118} for the proof):

\begin{lem*}
The differential of the volume density
\begin{equation*}
\left\{ \begin{array}{ccl}
\Ga(S^2_+T^*M) &\to &\Gamma(\vol(M)) \\
g &\mapsto &\on{vol}(g)
\end{array}\right.\end{equation*}
is given by
\begin{equation*}
D_{(g,m)} \on{vol}(g) = \frac{1}{2}\on{Tr}(g\i.m)\on{vol}(g).
\end{equation*}
\end{lem*}

\subsection{Metric on tensor fields}\label{nmb:2.4}

A metric on a space of tensor fields is defined by integrating the appropriate metric on the 
tensor space with respect to the volume density:
$$\wt g ^r_s(h,k)=\int_M g^r_s\big(h(x),k(x)\big)\on{vol}(g)(x)$$
for $h,k \in \Gamma(T^r_sM)$.
According to Section~\ref{nmb:2.1}, if $h$ and $k$ 
are tensor fields of type 
$\left(\begin{smallmatrix}0\\2\end{smallmatrix}\right)$ and $h$ or $k$ is symmetric, then
\begin{align*}
\wt g^0_2 (h,k) = \int_M \on{Tr}(g\i h(x) g\i k(x)) \on{vol}(g)(x).
\end{align*}

\subsection{Covariant derivative on $M$}\label{nmb:2.5}

We will use covariant derivatives on vector bundles as explained in \cite[especially Section~19.12]{MichorH}.
Let $X$ be a vector field on $M$. 
The Levi-Civita covariant derivative $\na_X$ on $(M, g)$ can be extended uniquely to an operator 
on the space $\Gamma(T^r_sM)$ of all tensor fields on $M$. This covariant derivative
depends on the metric $g$.

We define its derivative with respect to $g$ as
\begin{equation}
N^r_s(m)=N^r_s(g,m)=D_{(g,m)}\na,, 
\end{equation}
where 
$$\nabla \in L\big(\Ga(T^r_s M), \Ga(T^r_{s+1} M)\big)$$
and where $m$ is a tangent vector to  $\mathcal M(M)$ with foot point $g$.
The operator $N^r_s(m) \in \Ga\big(L(T^r_s M,T^r_{s+1}M)\big)$ is tensorial since 
$$D_{(g,m)} \na(fh) = D_{(g,m)} (df\otimes h+f \na h) = f D_{(g,m)} \na h$$
holds for $f \in C^\infty(M)$ and $h \in \Ga(T^r_s M)$.
In abstract index notation one has
\begin{equation}\label{eq:n10}
\big(N^1_0(m)\big)\vp^{\phantom{j}i}_{jk} = \frac12 g^{il} \big((\na m)_{jkl}+(\na m)_{kjl}-(\na m)_{ljk}\big),
\end{equation}
as can be seen from the formula \cite[theorem~1.174]{Besse2008}:
$$g\big(D_{(g,m)}(\na_X Y),Z\big) =\! 
\frac12 \big((\na_X m)(Y,Z)+(\na_Y m)(X,Z)-(\na_Z m)(X,Y)\big).$$
Furthermore, $(N^0_1(m))^i_{jk}=-(N^1_0(m))^i_{kj}$ since one has for $\al \in \Om^1(M)$ and $X,Y \in \X(M)$: 
\begin{align*}
\bigl(N^0_1(m)\al\bigr)(X,Y) &= (D_{(g,m)} \na_X \al)(Y) = D_{(g,m)} \big(d(\al(Y)).X-\al(\na_X Y) \big) 
\\&
= -\al(D_{(g,m)} \na_X Y) = -\bigl(N^1_0Y\bigr)(\al,X).
\end{align*}
Since $\na_X$ is a derivation on tensor products, one gets a similar property for $N^r_s(m)$: 
\begin{multline}\label{eq:npq}
\big(N^r_s(m)\big)\vp^{\vp}_j\vp^{i_1}_{k_1}\vp^{\ldots i_r}_{\ldots k_r}\vp^{i_{r+1}}_{k_{r+1}}
\vp^{\ldots i_{r+s}}_{\ldots k_{r+s}}
=\\=
\big(N^1_0(m)\big)\vp^{\phantom{j}i_1}_{jk_1} \de\vp^{i_2}_{k_2} \ldots \de\vp^{i_{r+s}}_{k_{r+s}}+\ldots+
\de\vp^{i_1}_{k_1} \ldots \de\vp^{i_{r+s-1}}_{k_{r+s-1}} \big(N^0_1(m)\big)\vp^{\phantom{j}i_{r+s}}_{jk_{r+s}},
\end{multline}
where one has $N^1_0$ in the first $r$ summands and $N^0_1$ in the last $s$ summands.


\subsection{The adjoint of the covariant derivative}\label{nmb:2.6}

The covariant derivative, seen as a mapping
$\na: \Gamma(T^r_sM) \to  \Gamma( T^r_{s+1}M)$
admits an adjoint 
$\na^*:\Gamma( T^r_{s+1}M)\to \Gamma(T^r_sM)$
with respect to the metric $\wt g$, i.e.: 
$\wt g^r_{s+1}(\na B, C)= \wt g^r_s(B, \na^* C).$
It is given by
$\na^*B=-\on{Tr}^g(\na B), $
where the trace contracts the first two tensor slots. 
This formula is proven in \cite{Michor119}.

\subsection{Second covariant derivative}\label{nmb:2.7}

When the covariant derivative is seen as a mapping
$\na: \Gamma(T^r_s M) \to \Gamma(T^r_{s+1}M),$
then the \emph{second covariant derivative} is simply 
$\na \na=\na^2: \Gamma(T^r_s M) \to \Gamma(T^r_{s+2}M).$
For $X,Y \in \X(M)$, it is given by
$\na^2_{X,Y} =\iota_Y \iota_X \na^2 =
\iota_Y \na_X \na =
\na_X\na_Y -\na_{\na_XY} .$
Higher covariant derivatives are defined accordingly.

\subsection{Laplacian}\label{nmb:2.8}

The \emph{Bochner-Laplacian} is defined as 
$\Delta h := \na^*\na h = - \on{Tr}^g(\na^2 h).$
It can act on all tensor fields $h$, and it respects the degree of the tensor field it is acting on. 
Using \ref{nmb:2.5} we get: 
\begin{lem*}
The differential of the Laplacian acting on
$\binom rs$-tensors is given by:
\begin{align*}
D_{(g,m)} \De h &= - D_{(g,m)} \Tr^g(\na^2 h) \\&=
\Tr(g\i m g\i \na^2 h) - \Tr^g\big(N^r_{s+1}(m) \na h\big) - \Tr^g\big(\na N^r_s(m) h\big).
\end{align*}
\end{lem*}
Here the trace contracts the first two tensor slots, for example
$$\Tr(g\i m g\i \na^2 h) = g^{ij} m_{jk} g^{kl}\nabla^2_{li}h.$$

\subsection{Curvature}\label{nmb:2.9}
The Riemann curvature tensor is given by $$R(X,Y)Z=\na_X\na_YZ-\na_Y\na_XZ-\na_{[X,Y]}Z.$$ 
The Ricci  tensor field $\Ric(X,Y)$ is the trace of $Z\mapsto R(Z,X)Y$. The scalar curvature is
$\Scal=\Tr^g(\Ric)$.

\begin{lem*} {\rm\cite[theorem 1.174]{Besse2008}} The differential of the scalar curvature 
\begin{equation*}
\left\{ \begin{array}{ccl}
\Ga(S^2_+T^*M) &\to &C^{\infty}(M),\\
g &\mapsto &\Scal
\end{array}\right.\end{equation*}
is given by
\begin{equation*}
D_{(g,m)} \Scal = \Delta(\Tr(g\i.m))+\na^*(\na^*(m))-g^0_2(\Ric,m).
\end{equation*}
\end{lem*}

\section{Riemannian metrics on the manifold of Riemannian metrics}\label{nmb:3}

Let $P_g:\Ga(S^2T^*M) \to \Ga(S^2T^*M)$
be a positive, symmetric, bijective pseudo-differential operator of order $2p$ 
depending smoothly on the metric $g$.
Then the operator $P$  induces a  metric on the manifold of Riemannian metrics, namely
$$G^P_g(h,k)=\int_M g^0_2(P_gh,k) \vol(g)=\int_M \Tr(g\i. P_gh. g\i. k) \vol(g).$$ 

\subsection{Geodesic equation}\label{nmb:3.1}

Given $(1,2)$-tensors $H$ and $K$ on $\mathcal M(M)$ such that
$$D_{(g,m)}G_g^P(h,k)=G_g^P(K_g(h,m),k)=G_g^P(m,H_g(h,k)),$$
the geodesic equation is given by the following variant of the Christoffel symbols 
$$g_{tt}=\frac12 H_g(g_t,g_t)-K_g(g_t,g_t),$$
see \cite{Michor107,Michor118,Michor119}.

We will now compute the metric gradients $H$ and $K$. 
The calculations at the same time show the existence of the gradients. 
For this aim, let  $m,h,k \in T_g\Met$ be constant vector fields on $\mathcal M(M)$. Using the formula for the variation of the volume density 
from Section \ref{nmb:2.3} we get
\begin{align*}
G_g^P&(K_g(h,m),k)=D_{(g,m)}G_g^P(h,k) =
\\&
=D_{(g,m)}\int_M \Tr(g\i.Ph.g\i.k) \vol(g)
\\&
=\int_M \Tr\big((D_{(g,m)}g\i).Ph.g\i.k\big) \vol(g)
\\&\quad
+\int_M \Tr\big(g\i.(D_{(g,m)}P)h.g\i.k\big) \vol(g)
\\&\quad
+\int_M \Tr\big(g\i.Ph.(D_{(g,m)}g\i).k\big) \vol(g)
\\&\quad
+\int_M \Tr\big(g\i.Ph.g\i.k\big)D_{(g,m)}\vol(g)
\\&
=\int_M\Big[- \Tr\big(g\i.m.g\i.Ph.g\i.k\big) 
\\&\qquad\qquad
+ \Tr\big(g\i.(D_{(g,m)}P)h.g\i.k\big) 
-\Tr\big(g\i.Ph.g\i.m.g\i.k\big) 
\\&\qquad\qquad\qquad\qquad\qquad\quad
+\Tr\big(g\i.Ph.g\i.k\big)\frac12\Tr(g\i.m)\Big]\vol(g)
\\&
=\int_M g^0_2\Big(-m.g\i.Ph+(D_{(g,m)}P)h-Ph.g\i.m
\\&\qquad\qquad\qquad\qquad\qquad\qquad\qquad\qquad
+\frac12\Tr(g\i.m).Ph,k\Big) \vol(g).
\end{align*}
Therefore the $K$-gradient is given by 
\begin{multline*}
K_g(h,m)=
P\i\Big[-m.g\i.Ph\, +
\\
+(D_{(g,m)}P)h-Ph.g\i.m+\frac12\Tr(g\i.m).Ph\Big].
\end{multline*}
To calculate the $H$-gradient we will assume that there exists an \emph{adjoint} 
in the following sense
\begin{equation*}
\boxed{
\int_M g^0_2\big((D_{(g,m)}P)h,k\big) \vol(g)=\int_M g^0_2\big(m,(D_{(g,.)}Ph)^*(k)\big) \vol(g)
}\tag{1}\end{equation*}
which is smooth in $(g,h,k)$ and bilinear in $(h,k)$.
The existence of the adjoint needs to be checked for each specific operator $P$, 
usually by partial integration.  
Using the adjoint we can rewrite the equation above as follows:
\begin{align*}
&G_g^P(H_g(h,k),m)=(D_{(g,m)}G_g^P)(h,k)=D_{(g,m)}\int_M g^0_2(Ph,k) \vol(g)
\\&\quad
=\int_M g^0_2\Big(-m.g\i.Ph+(D_{(g,m)}P)h-Ph.g\i.m
\\&\qquad\qquad\qquad\qquad\qquad\qquad
+\frac12\Tr(g\i.m).Ph,k\Big) \vol(g)
\\&\quad
=\int_M \Big(g^0_2\big(m,-Ph.g\i.k\big)+g^0_2\big(m,(D_{(g,.)}Ph)^*(k)\big)
\\&\qquad\qquad
+g^0_2\big(m,-k.g\i.Ph\big)
+\frac12g^0_2\big(m,g.\Tr(g\i.Ph.g\i.k)\big)\Big)\vol(g)
\end{align*} 
Here we can easily read off the $H$-gradient: 
\begin{multline*}
H_g(h,k)=
P\i\Big[(D_{(g,.)}Ph)^*(k) -
\\
-Ph.g\i.k-k.g\i.Ph+\frac12.g.\Tr(g\i.Ph.g\i.k)\Big].
\end{multline*}
Therefore the geodesic equation on the manifold of Riemannian metrics reads as:
$$\boxed{\begin{aligned}
g_{tt}&=\frac12 H_g(g_t,g_t)-K_g(g_t,g_t)
\\&
=P\i\Big[\frac12 (D_{(g,.)}Pg_t)^*(g_t)+\frac14.g.\Tr(g\i.Pg_t.g\i.g_t)
+\frac12 g_t.g\i.Pg_t
\;\\&\qquad\qquad
+\frac12 Pg_t.g\i.g_t-(D_{(g,g_t)}P)g_t-\frac12\Tr(g\i.g_t).Pg_t\Big]
\end{aligned}}$$
We can rewrite this equation to get it in a slightly more compact form:
\begin{equation*}
\boxed{\begin{aligned}
(Pg_{t})_t&=(D_{(g,g_t)}P)g_t+Pg_{tt}
\\&=
\frac12(D_{(g,.)}Pg_t)^*(g_t)+\frac14.g.\Tr(g\i.Pg_t.g\i.g_t)
\\&\qquad\quad
+ \frac12 g_t.g\i.Pg_t+\frac12 Pg_t.g\i.g_t-\frac12\Tr(g\i.g_t).Pg_t\;
\end{aligned}
\tag{2}}\end{equation*}

\subsection{Well-posedness of some geodesic equations}\label{nmb:3.2}
For any fixed background Riemann metric $\hat  g$ on $M$ and 
its Levi-Civita covariant derivative $\hat  \na$, 
the {\it Sobolev space} $H^k(S^2T^*M)$ is the Hilbert space completion of the space 
$\Ga(S^2T^*M)$ of smooth sections, in the Sobolev norm
$$
\|h\|_k^2 = \sum_{j=0}^k \int_M  \hat  g^0_{2+j} ((\hat  \na)^j h, (\hat  \na)^j h)\vol(\hat  g).
$$
The topology of the Sobolev space does not depend on the choice of $\hat  g$; the resulting norms are
equivalent. 
See \cite{Shubin1987}  for more information.
The following results hold:
\begin{itemize}
       \item {\it Sobolev lemma.} If $k>\frac{\dim(M)}{2}$ then the identity on $\Ga(S^2T^*M)$ extends to an
       injective bounded linear mapping $H^{k+p}(S^2T^*M)\to C^p(S^2T^*M)$ where $C^p(S^2T^*M)$ carries the 
       supremum norm of all derivatives up to order $p$.
       \item {\it Module property of Sobolev spaces.}
	If $k>\frac{\dim(M)}{2}$ then the evaluation $H^k(L(S^2T^*M,S^2T^*M))\x H^k(S^2T^*M)\to H^k(S^2T^*M)$ is
       bounded and bilinear. Likewise all other point wise contraction operations are multilinear bounded
       operations. See \cite{EichhornFricke1998}, or \cite[1.3.12]{Eichhorn2007}.
\end{itemize}
The Sobolev lemma allows us to define the Sobolev space $\Met^k(M):=H^k(S^2_+T^*M)$
for $k>\frac{\dim(M)}{2}$. 

\subsection*{Assumptions} {\it
In the following we assume the natural condition that $h\mapsto P_gh$ is an elliptic and self-adjoint 
pseudo-differential operator of order $2p\geq 0$. Then 
it is Fredholm and it has vanishing index by \cite[theorem 26.2]{Shubin1987}.
Thus it is invertible and $g\mapsto P_g\i$ is a smooth mapping
$$H^k(S^2_+T^*M)\to L(H^{k}(S^2T^*M),H^{k+2p}(S^2T^*M))$$ by the implicit function theorem on Banach 
spaces.

We assume that $(D_{(g,.)}Ph)^*(m)$ exists and is a linear pseudo-differen\-tial operator of order $2p$ in $m,h$.

As (non-linear) mappings in the foot point $g$, we assume that 
$P_g h$, $(P_g)\i h$, $(D_{(g,.)}Ph)^*(m)$ are compositions of operators of the following type: 
\begin{enumerate}[(a)]
\item Non-linear differential operators of order $l\leq 2p$, i.e.
$$A(g)(x)=A\big(x,g(x),(\hat\na g)(x),\ldots,(\hat\na^l g)(x)\big),$$
\item Linear pseudo-differential operators of order $\leq 2p$,
\end{enumerate}
such that the total (top) order of the composition is $\leq 2p$.
}


\begin{thm*}
Let the assumptions above hold.
Then for $k>\frac{\dim(M)}{2}+1$, the initial value problem for the geodesic equation 
\thetag{\ref{nmb:3.1}.2}
has unique local solutions in the Sobolev manifold $\Met^{k+2p}(M)$ of
$H^{k+2p}$-metrics. The solutions depend $C^\infty$ on $t$ and on the initial
conditions $g(0,\;.\;)\in \Met^{k+2p}(M)$ and $g_t(0,\;.\;)\in H^{k+2p}(S^2T^*M)$.
The domain of existence (in $t$) is uniform in $k$ and thus this
also holds in $\Met(M)$.

Moreover, in each Sobolev completion $\Met^{k+2p}(M)$, the Riemannian exponential mapping $\exp^{P}$ exists 
and is smooth on a neighborhood of the zero section in the tangent bundle, 
and $(\pi,\exp^{P})$ is a diffeomorphism from a (smaller) neighborhood of the zero 
section to a neighborhood of the diagonal in 
$\Met^{k+2p}(M)\x \Met^{k+2p}(M)$. 
All these neighborhoods are uniform in $k>\frac{\dim(M)}{2}$
and can be chosen $H^{k_0+2p}$-open, where $k_0>\frac{\dim(M)}{2}$.
Thus all properties of the exponential 
mapping continue to hold in $\Met(M)$.
\end{thm*}

This proof is an adaptation of \cite[section 4.2]{Michor119}.

\begin{proof}
We consider the geodesic equation as the flow equation of a smooth
($C^\infty$) vector field $X$ on the open set 
$$
\Met^{k+2p}\x H^k(S^2T^*M) \subset H^{k+2p}(S^2T^*M)\x H^k(S^2T^*M).
$$

We now write the geodesic equation as the flow equation 
of an autonomous smooth vector field $X=(X_1,X_2)$ on $\Met^{k+2p}\x H^k$, as 
follows (using 
\thetag{\ref{nmb:3.1}.2}:
\begin{align*}
g_t &= (P_g)\i h =:X_1(g,h)
\\
h_t&= 
\frac12\big((D_{(g,.)}P_g) (P_g)\i h\big)^*((P_g)\i h)+\frac14.g.\Tr(g\i.h.g\i.(P_g)\i h)
\\&\quad
+ \frac12 (P_g)\i h.g\i.h+\frac12 h.g\i.(P_g)\i h-\frac12\Tr(g\i.(P_g)\i h).h
\tag{1}\\&
=: X_2(g,h)
\end{align*}
For $(g,h)\in \Met^{k+2p}\x H^k$ we have $(P_g)\i h\in H^{k+2p}$.
Thus 
a term by term investigation of \thetag{1}, using the assumptions on the orders, shows that 
$X_2(g,h)$
is smooth in  $(g,h)\in \Met^{k+2p}\x H^k$ with values in $H^{k}$.
Likewise $X_1(g,h)$ is smooth in $(f,h)\in \Met^{k+2p}\x H^k$ with values in 
$H^{k+2p}$. Thus by the theory of smooth ODE's on Banach spaces,
the flow $\on{Fl}^k$ exists on $\Met^{k+2p}\x H^k$ and is smooth in $t$ and the initial 
conditions for fixed $k>\frac{\dim(M)}{2}+1$.

We choose $C^\infty$ initial conditions $g_0=g(0,\quad)$ and
$h_0=P_{g_0} g_t(0,\quad)=h(0,\quad)$ for the flow equation \thetag{1} in 
$\Met(M)\x \Ga(S^2T^*M)$. Suppose the 
trajectory $\on{Fl}^k_t(g_0,h_0)$ of
$X$ through these initial conditions in $\Met^{k+2p}\x H^k$ maximally exists for $t\in (-a_k,b_k)$, 
and the trajectory $\on{Fl}^{k+1}_t(g_0,h_0)$ in $\Met^{k+1+2p}\x H^{k+1}$ 
maximally exists for $t\in(-a_{k+1},b_{k+1})$ with $a_{k+1}<a_k$ and $b_{k+1}<b_k$, say. Since
solutions are unique, $\on{Fl}^{k+1}_t(g_0,h_0)=\on{Fl}^{k}_t(g_0,h_0)$ for
$t\in (-a_{k+1,}b_{k+1})$. 
We now apply the background derivative $\hat  \na$ to both equations \thetag{1}:
\begin{align*}
(\hat  \na g)_t &= \hat  \na g_t =  \hat  \na X_1(g,h) 
\\
(\hat  \na h)_t &=\hat  \na h_t =  \hat  \na X_2(g,h) 
\end{align*}
We claim that for  $i=1,2$ we have 
$$
\hat  \na X_i(g,h) 
= X_{i,1}(g,h)(\hat  \na^{2p+1}g) + X_{i,2}(g,h)(\hat  \na^{2p+1}h) + X_{i,3}(g,h)
$$
where all $X_{i,j}(g,h)(l)$ and $X_{i,3}(g,h)$ ($i,j=1,2$) are smooth in all variables, 
of highest order $2p$ in $g$ and $h$, 
linear and algebraic (i.e., of order 0) in $l$.
This claim follows from the assumptions: (b) For a linear pseudo differential operator $B$ of 
order $q$ 
the commutator $[\na_Y,B]$ is a pseudo differential operator of order $q$ again for any vector field $Y$. 
(a) For a local operator we can apply the chain rule: The derivative of order $2p+1$ of $g$ appears only linearly.

We write 
$\hat  \na^{2p+1}g = \hat  \na^{2p}\tilde g$ and 
$\hat \na^{2p+1}h = \hat \na^{2p}\tilde h$ for the highest derivatives only.
Then $\tilde g$ and $\tilde h$ satisfy
\begin{align*}
\tilde g_t &= X_{1,1}(g,h)(\hat \na^{2p}\tilde g) + 
X_{1,2}(g,h)(\hat \na^{2p}\tilde h) + X_{1,3}(g,h)
\\
\tilde h_t &= X_{2,1}(g,h)(\hat \na^{2p}\tilde g) + 
X_{2,2}(g,h)(\hat \na^{2p}\tilde h) + X_{2,3}(g,h).
\end{align*}
This ODE is inhomogeneous bounded and linear
in $(\tilde g,\tilde h)\in \Met^{k+2p}\x H^k$ with
coefficients bounded linear operators on $H^{k+2p}$ and $H^k$, respectively.  
These coefficients are $C^\infty$ functions of $(g,h) \in \Met^{k+2p}\x H^k\subset C^1$ 
which we already know on the interval $(-a_k,b_k)$. 
This equation therefore has a solution $(\tilde g(t,\quad),\tilde h(t,\quad)) \in \Met^{k+2p}\x H^k$ for all $t$ for which the
coefficients exists, thus for all $t\in (-a_k,b_k)$. 
Obviously, $(\tilde g,\tilde h)=(\hat \na g,\hat \na h)$ for $t \in (-a_{k+1},b_{k+1})$.
By continuity this holds also for $t \in [-a_{k+1},b_{k+1}]$ which contradicts that the interval $(-a_{k+1},b_{k+1})$ is maximal. 
We can iterate this and conclude that the flow of $X$ exists in
$\bigcap_{m\ge k} \Met^{m+2p}\x H^{m}= \Met\x \Ga$. 

It remains to check the properties of the Riemannian exponential mapping $\exp^P$.
It is given by $\exp^P_{g}(h)= c(1)$ where $c(t)$ is the geodesic emanating from  
value $g$ with initial velocity $h$. 
From the form  $g_{tt}=\frac12 H_g(g_t,g_t)-K_g(g_t,g_t) =: \Ga_g(g_t,g_t)$  
(see subsection \ref{nmb:3.1}), namely linearity in $g_{tt}$ and 
bilinearity in $g_t$, and from local existence and uniqueness on each space $\Met^{k+2p}(M)$ 
the properties claimed follow:
see for example \cite[22.6 and 22.7]{MichorH} for a detailed proof in terms of the spray
vector field $S(g,h)=(g,h;h,\Ga_g(h,h))$ on a finite dimensional manifold. This proof carries over to infinite dimensional 
convenient manifolds without any change in notation. 
So we check this on the largest of this spaces $\Met^{k_0}(M)$ (with the smallest $k$). 
Since the spray on $\Met^{k_0}(M)$ restricts to the spray on each $\Met^{k+2p}(M)$,
the exponential mapping $\exp^P$ and the inverse $(\pi,\exp^P)\i$ on $\Met^{k_0}(M)$ restrict 
to the corresponding mappings on each $\Met^{k+2p}(M)$. Thus the neighborhoods of existence are 
uniform in $k$.
\end{proof}

\subsection{Conserved Quantities}\label{nmb:3.3}
Consider the right action of the diffeomorphism group $\on{Diff}(M)$ on $\Met(M)$ given by
$(g,\phi) \mapsto \phi^*g$ with fundamental vector field 
$$\zeta_X(g)=\L_Xg=2\on{Sym}\na(g(X)).$$
For a proof of the last equality see \cite[section 1]{Besse2008}.
If the metric $G^P$ is invariant under this action, we have the following conserved quantities 
(see for example \cite{Michor118}):
\begin{align*}
\on{const}&=G^P(g_t,\zeta_X(g))=\int_M g^0_2\big(Pg_t,2\on{Sym}\na(g(X))\big)\vol(g)\\&
=2\int_M g^0_1\big(\na^*\on{Sym}Pg_t,g(X)\big)\vol(g)=2\int_M (\na^*Pg_t)(X)\vol(g)\\&
=2\int_M g\big(g\i\na^*Pg_t,X\big)\vol(g)
\end{align*}
Since this equation holds for all vector fields $X$ this yields
$$\boxed{(\na^*Pg_t)\on{vol}(g)\in \Ga(T^*M\otimes_M\on{vol}(M))\text{ is const.\ in time}.}$$
The geometric interpretation of this conserved quantity is carried by the expression 
$G^P(g_t,\zeta_X)$. 
After normalization this gives a formula for the cosine of the angle between the geodesic and any
vector field $\zeta_X$. If the constant vanishes  
then this geodesic is $G^P$-perpendicular to each $\on{Diff}(M)$-orbit it meets.


\subsection{Geodesics of pure scalings}\label{nmb:3.4}
In this section we want to investigate when $r(t)g_0$ is a geodesic for some real function $r$ and 
some fixed metric $g_0$. This will help us to determine the geodesic completeness of the space $\Met(M)$
under various metrics. 

\begin{lem*}
Let $g_0 \in \Met(M)$ and $\N=\NN=\{r g_0: r >0\} \subset \Met(M)$.
If $P$ viewed as a $\binom{1}{1}$-tensor field on $\Met(M)$ `restricts' to 
the submanifold $\N$ in the sense that $P_gh$ is tangential to $\N$ for all 
$g \in \N$ and $h \in T_g\N$,
then the following statements are equivalent.
\begin{enumerate}[(a)]
	\item $\N$ is totally geodesic.
	\item $(D_{(g,\cdot)}P h)^*(k)$ is tangential to $\N$ for all $g \in \N$ and $h,k \in T_g\N$.
	\item $(D_{(g,m)}P)(h)$ is $\wt g^0_2$-normal to $\N$ for all $g \in \N$, $h \in T_g\N$, 
		$m \in T_g\Met(M)$ such that $m$ is $\wt g^0_2$-normal to $\N$.
\end{enumerate}
If $P$ restricts to $\N$ and (a)-(c) hold, then there are $\Psi,f:\R_{>0} \to \R$ such that
$$P_{r g_0}(g_0)= \Psi(r)g_0, \qquad ((D_{(r g_0,\cdot)}P) g_0)^*(g_0)=f(r)g_0$$ holds for all $r>0$. 
Then the path $g(t, \cdot)= r(t).g_0$ is a geodesic in $\Met(M)$ if and only if the function $r$ satisfies 
\begin{align*}
	r'' \Psi(r) = r'^2\Big(\frac12 f(r)-\Psi'(r)  +\big(1-\dim(M)/4\big) \Psi(r) r\i \Big) \enspace .
\end{align*}
Along these geodesics the conserved quantity vanishes, i.e.,
$$
(\nabla^*Pg_t)\vol(g)=0\enspace.
$$
\end{lem*}
\begin{rem*}
Note that $(D_{(g,m)}P)(h)$ and $(D_{(g,\cdot)}P h)^*(k)$ are tensorial in $h,k$ and that
for $g\in\N$, all tangent vectors in $T_g\N$ can be written as real multiples of $g$. 
Therefore conditions ($b$) and ($c$) are equivalent to 
\begin{enumerate}[(a')]\setcounter{enumi}{1}
 	\item $(D_{(g,\cdot)}P g)^*(g)$ is tangential to $\N$ for all $g \in \N$.
	\item $(D_{(g,m)}P)(g)$ is $\wt g^0_2$-normal to $\N$ for all $g \in \N$ and 
		$m \in T_g\Met(M)$ such that $m$ is $\wt g^0_2$-normal to $\N$.
\end{enumerate}
\end{rem*}

\begin{proof}
The submanifold $\N$ is totally geodesic iff  
$\frac12 H_g(h,k)-K_g(h,k)$ is tangential to $\N$
for all $g \in \N$ and $h,k \in T_g\N$.
We now look at the formulas for $H$ and $K$ from Section~\ref{nmb:3.1}. 
Since $P_g$ is bijective and preserves $T_g\N$, the above condition is equivalent to 
$P_g(\frac12 H_g(h,k)-K_g(h,k))$ being tangential. 
A term-by-term investigation shows that this is the case if and only if 
$((D_{(g,\cdot)}P)h)^*(k)$ is tangential, in which case it can be expressed using a function $f$. 
A test for the latter condition is $$\wt g^0_2\big(((D_{(g,\cdot)}P)h)^*(k),m\big)=\wt g^0_2\big((D_{(g,m)}P)h,k\big)=0$$ for all 
$m\in T_g\Met(M)$ that are $\wt g^0_2$-normal to $\N$.  Equivalently, 
$(D_{(g,m)}P)h$ has to be $\wt g^0_2$-normal to $\N$ 
whenever $m$ is $\wt g^0_2$-normal to $\N$ and $h$ is tangential to $\N$. 

It remains to check  the form of the geodesic equation.
We use the geodesic equation (2) from Section~\ref{nmb:3.1} and substitute
$$g=r(t)g_0, \quad g_t=r'(t)g_0, \quad P_g(g_t)=\Psi(r(t))r'(t)g_0.$$ 
Dropping the dependence on $t$ in our notation we 
get for the left hand side of the geodesic equation:
\begin{align*}
\p_t(P_g g_t) &= r'' \Psi(r) g_0+ \Psi'(r) r'^2 g_0
\end{align*}
The previous substitutions and 
$$((D_{(g,.)}P)g_t)^*(g_t)=f(r(t)) r'(t)^2 g_0$$
yield the right-hand side of the geodesic equation:
\begin{align*}
	&\frac12(D_{(g,.)}Pg_t)^*(g_t)+\frac14.g.\Tr(g\i.Pg_t.g\i.g_t)+ \frac12 g_t.g\i.Pg_t\\
	&\qquad\qquad\qquad+\frac12 Pg_t.g\i.g_t-\frac12\Tr(g\i.g_t).Pg_t\\
	&\qquad=\frac12 f(r) r'^2 g_0+\frac14 \Psi(r) \dim(M) r\i r'^2  g_0+ \frac12 \Psi(r) r\i r'^2  g_0\\
	&\qquad\qquad\qquad +\frac12 \Psi(r) r\i r'^2 g_0-\frac12 \Psi(r) \dim(M) r\i r'^2 g_0\\
	&\qquad=\frac12 f(r) r'^2 g_0+\big(1-\dim(M)/4\big) \Psi(r) r\i r'^2  g_0. 
\end{align*}
For the conserved quantity we calculate:
\begin{align*}
(\nabla^*Pg_t)\vol(g)&= \Tr\big(g\i\nabla^{g}(r'(t)g_0)\big)\vol(g)
\\&
= \frac{r'(t)}{r(t)}\Tr\big(g\i\nabla^{g}(g)\big)\vol(g)=0\;. 
\end{align*}
\end{proof}

\subsection{Length of pure scalings}\label{nmb:3.5}
\begin{lem*}
Given $g_0$ such that $P_{rg_0}(g_0)=\Psi(r).g_0$ the length of the curve $g(r)=r g_0$  for $r\in[0,1]$ is given by
\begin{align*}
\on{Len}_0^1(g) &= 
\sqrt{\dim(M)\Vol(g_0)} \int_0^1 \sqrt{\Psi(r) r^{\dim(M)/2-2}}dr \enspace .
\end{align*}
If $\Psi(r)=O(r^{\al})$ 
for some $\al>-\dim(M)/2$, then 
$\NN \subset \Met(M)$ is an incomplete metric space under $G^P$.
If in addition $P$ and $g_0$ satisfy the conditions of Lemma~\ref{nmb:3.4}, then 
$(\Met(M), G^P)$ is geodesically incomplete.
\end{lem*}

Note that $(\Met(M), G^P)$ is always an incomplete metric space since it does not contain Sobolev 
class $H^p$ metrics.

\begin{proof}
For the length of the curve we calculate:
\begin{align*}
\on{Len}_0^1(g) &= \int_0^1 G^P_{r.g_0}(g_0,g_0)^{1/2}\,dr 
\\&
= \int_0^1 \Big(\int_M \Tr((rg_0)\i.P_{rg_0}(g_0).(rg_0)\i.g_0) \vol(rg_0)\Big)^{1/2}dr
\\&
= \int_0^1 r^{\dim(M)/4-1}\Big(\int_M \Tr((g_0)\i.P_{rg_0}(g_0)) \vol(g_0)\Big)^{1/2}dr.
\end{align*}
Using  the  assumption $P_{rg_0}(g_0)=\Psi(r).g_0$,
we can compute this as 
\begin{align*}
\on{Len}_0^1(g) &= 
\int_0^1 r^{\dim(M)/4-1}\sqrt{\dim(M)}\Big(\int_M \Psi(r) \vol(g_0)\Big)^{1/2}dr\enspace.
\end{align*}
Note that the metric space $(\Met(M), G^P)$ is geodesically incomplete 
if $\NN$ contains a geodesic in $\Met(M)$ which connects $g_0$ to 0 in finite time. 
\end{proof}

\section{Special cases of $P$}\label{nmb:4}
In this section we present various interesting examples of metrics. These special choices are motivated 
by related  metrics on spaces of immersions and shape spaces, see 
\cite{Michor119,Michor118,MBMB2011}. We will use the notation $n=\dim(M)$ for all of this section.

\subsection{The $H^0$-metric}\label{nmb:4.1}
The simplest and most natural example is the operator $P$ of order zero given by
$P_g(h)=h$ for $g \in \Met(M)$ and $h \in T_g\Met(M)$.
With this choice of $P$, the metric $G^P$ equals $\wt g^0_2$.
It is the so called $L^2$-metric or $H^0$-metric, 
which is well studied as mentioned in the introduction. 
We can easily read off the geodesic equation from the previous section:
$$\boxed{g_{tt}=\frac14.g.g^0_2(g_t,g_t)+g_t.g\i.g_t-\frac12\Tr(g\i.g_t).g_t.}$$
This coincides with the equation derived in \cite{FreedGroisser89} and \cite{Gil-MedranoMichor91}.
All conditions from \ref{nmb:3.2} are obviously satisfied. Thus the geodesic equation is well-posed.
Here the geodesic equation evolves in each set $S^2_+ T_x^*M$ separately. 
The conserved quantities have the form
$$\boxed{(\nabla^* g_t)\on{vol}(g)\in \Ga(T^*M\otimes_M\on{vol}(M))\text{ is const.\ in time}.}$$

The conditions of Lemma~\ref{nmb:3.4} are obviously satisfied for all 
$g_0\in\Met(M)$ and we get again the result from \cite{Gil-MedranoMichor91} that 
$\NN$ is the image of a geodesic. 
The geodesic is $r(t).g_0$ where $r(t)$ satisfies 
\begin{align*}
r''(t)&=\frac{r'(t)^2}{r(t)} \big(1-\frac{n}{4} \big)\;,\text{  i.e., }\;
r(t)= \Big(t (r(1)^{n/4}-r(0)^{n/4})+r(0)^{n/4} \Big)^{4/n}  .
\end{align*}
This geodesic connects $g_0$ with 0 in finite time. Thus
it follows that the space $(\Met(M),\wt g^0_2)$ is geodesically incomplete.

\subsection{Conformal metrics}\label{nmb:4.2}
Here we consider metrics of the form
$$G^{P}_g(h,k)=\Phi(\Vol(g))\int_M g^0_2(h,k) \vol(g),$$ 
where $\Phi\in C^\infty(\mathbb R_{>0},\mathbb R_{>0})$ and $\on{Vol}(g)=\int_M\on{vol}(g)$.
To calculate the adjoint we will use the variational formula for the volume form from section~\ref{nmb:2.3}:
\begin{align*}
&\int_M g^0_2\big(m,(D_{(g,.)}Ph)^*(k)\big) \vol(g)=\int_M g^0_2\big((D_{(g,m)}P)h,k\big) \vol(g)\\&\qquad
=\Phi'.(D_{(g,m)}\Vol).\int_M g^0_2\big(h,k\big) \vol(g)\\&\qquad
=\frac12\Phi'.\int_M \Tr(g\i.m)\vol(g).\int_M   g^0_2\big(h,k\big)\vol(g)\\&\qquad
=\frac12\int_M \Tr\Big(g\i.m.\Phi'.\int_M g^0_2(h,k) \vol(g)\Big)\vol(g)\\&\qquad
=\frac12 \int_M g^0_2\Big(m,\Phi'.g.\int_M g^0_2(h,k) \vol(g)\Big)\vol(g)
\end{align*}
Using this formula for the adjoint, the geodesic equation reads as:
$$\boxed{\begin{aligned}
&g_{tt}=\frac{\Phi'}{4\Phi}.g.\int_M g^0_2(g_t,g_t)\vol(g)+\frac{1}{4}.g.g^0_2(g_t,g_t)
+ g_t.g\i.g_t
\\&\qquad\qquad\qquad
-\frac{\Phi'}{2\Phi}.g_t.\int_M g_2^0(g_t,g) \vol(g)-\frac{1}{2}g^0_2(g_t,g).g_t\end{aligned}}$$
or
$$\boxed{\begin{aligned}
&(\Phi.g_{t})_t=\frac{\Phi'}{4}.g.\int_M g^0_2(g_t,g_t)\vol(g)+\frac{\Phi}{4}.g. g^0_2(g_t,g_t)\\
&\qquad\qquad\qquad+ \Phi. g_t.g\i.g_t-\frac{\Phi}{2} g^0_2(g_t,g).g_t\end{aligned}}$$
All conditions of theorem \ref{nmb:3.2} are satisfied. Thus the geodesic equation is 
well-posed and the geodesic exponential mapping exists and is a local diffeomorphism. 
Since the total volume $\Vol(M)$ does not depend on the point $x\in M$, the 
conserved quantities are: 
$$\boxed{\Phi(\Vol)\Tr(g\i \nabla g_t)\on{vol}(g)\in \Ga(T^*M\otimes_M\on{vol}(M))\text{ is const.\ in time}.}$$

Now we want to study again whether there exist metrics $g_0$ and  positive real function $r$ such that $r(t)g_0$ is a geodesic. 
Therefore we check whether the conditions of  Lemma~\ref{nmb:3.4} are satisfied.
$P$ obviously restricts to the submanifold $\NN$ for every $g_0\in\Met(M)$. 
Using again the variational formula for $\Vol$, we get   
\begin{multline*}
 \wt g^0_2\big((D_{(g,m)}P)g,g\big)= \frac12\Phi'.\int_M \Tr(g\i.m)\vol(g)\;\int_M   g^0_2\big(g,g\big)\vol(g)
 \\
= \frac12\Phi'\;\int_M \Tr(g\i.m.g\i.g)\vol(g)\; \int_M  n \vol(g) 
\\
= \frac{n}2\Phi'\Vol\; \wt g^0_2(m,g)=0\;,
\end{multline*}
if $m$ is $\wt g^0_2$-normal to $\NN$. Thus $\NN$ is a totally geodesic submanifold for any $g_0\in\Met(M)$.
For the corresponding functions $\Psi$ and $f$ we obtain:
\begin{align*}
P_{rg_0}g_0=\Psi(r)g_0&\text{ with } \Psi(r):= \Phi(r^{\tfrac{n}2} \Vol(g_0)) \,.\\
(\!(D_{(r g_0,\cdot)}P) g_0)^*(g_0)=f(r)g_0 &\text{ with }
f(r):=\frac{n}{2} \Phi'\big(r^{\tfrac{n}2}\Vol(g_0)\!\big)r^{\tfrac{n}{2}-1}\!\Vol(g_0).
\end{align*}
 The geodesic equation on  $\NN$ is then given by 
$$
r''\Phi(\Vol(rg_0))=\frac{r'^2}{r}\Big(-\frac{n}4 \Phi'(\Vol(r g_0))\Vol(r g_0)+\big(1-\frac{n}{4} \big) \Phi(\Vol(rg_0)) \Big).
$$

Let us now consider the special case $\Phi(\Vol)=\Vol^k$ for real $k$. Then the ODE for $r(t)$ simplifies to
$$
r''=\frac{r'^2}{r}\Big(1-\frac{n}{4}(k+1)\Big).
$$
with solution
$$r(t) = \Big(t (r(1)^{a}-r(0)^{a})+r(0)^{a} \Big)^{\frac{1}{a}}\quad \text{where}\quad a = \frac{n}{4}\big(1+k\big)   .$$
This geodesic connects $g_0$ with 0 in finite time if and only if $k>-1$. Thus
$(\Met(M),G^{\Ph(\Vol)})$ is geodesically incomplete if $\Ph(r)=O(r^k)$ for $r\searrow 0$, for some 
$k>-1$. Note that this would also follow from  Lemma~\ref{nmb:3.5}, since 
$\Psi(r)=\Phi\bigl(r^{\frac{n}{2}}\Vol(g_0)\bigr)$.

\subsection{Curvature weighted metrics}\label{nmb:4.3}
We consider metrics weighted by scalar curvature, 
$$G^{P}_g(h,k)=\int_M \Phi(\Scal(g)).g^0_2(h,k) \vol(g),$$
where $\Ph\in C^\infty(\mathbb R,\mathbb R_{>0})$.
Using the variational formula from section~\ref{nmb:2.9} we can calculate the adjoint as follows:
\begin{align*}
&\int_M g^0_2\big(m,(D_{(g,.)}Ph)^*(k)\big) \vol(g)=\int_M g^0_2\big((D_{(g,m)}P)h,k\big) \vol(g)
\\&
=\int_M \Phi'.(D_{(g,m)}\Scal) g^0_2\big(h,k\big) \vol(g)
\\&
=\int_M \Phi'.\Big(\Delta(\Tr(g\i.m))+\na^*(\na^*(m))-g^0_2(\Ric,m)\Big) g^0_2\big(h,k\big) \vol(g)
\\&  
=\int_M \Phi'.\Big[g^0_1\Big(\na\Tr(g\i.m),\na g^0_2(h,k)\Big)+g^0_1\Big(\na^*(m),\na g^0_2(h,k)\Big)
\\&\qquad\qquad\qquad\qquad\qquad\qquad
-g^0_2\Big(g^0_2(h,k)\Ric,m\Big) \Big]\vol(g)
\\&  
=\int_M \Phi'.\Big[\Tr(g\i.m).\na^*\na g^0_2(h,k)+g^0_2\Big(m,\na^2 g^0_2(h,k)\Big)
\\&\qquad\qquad\qquad\qquad\qquad\qquad
-g^0_2\Big(g^0_2(h,k)\Ric,m\Big) \Big]\vol(g)
\\&  
=\int_M \Phi'.g^0_2\Big(m,g.\De g^0_2(h,k)+\na^2 g^0_2(h,k)-g^0_2(h,k)\Ric\Big)\vol(g)
\end{align*}
Using the formula for the geodesic equation from section~\ref{nmb:3.1} yields
$$\boxed{\begin{aligned}
&(\Phi.g_{t})_t=\frac{\Phi'}{2}\Big(g.\De^g g^0_2(g_t,g_t)+\na^2 g^0_2(g_t,g_t)-g^0_2(g_t,g_t)\Ric\Big)\\&\qquad\qquad\qquad
+\frac{\Phi}{4}.g. g^0_2(g_t,g_t)+ \Phi. g_t.g\i.g_t-\frac{\Phi}{2} g^0_2(g_t,g).g_t.\end{aligned}}$$
The conditions of theorem~\ref{nmb:3.2} are violated and therefore it is not applicable. 
We do not know whether the geodesic equation is well-posed. The conserved quantities are given by
$$\boxed{\begin{aligned}
 &\nabla^* (\Phi(\Scal) g_t) \on{vol}(g)\\&= \bigg(\Phi'(\Scal)\Tr\big(g\i \on{dScal} \otimes g_t\big)+\Phi(\Scal)\Tr\big(g\i  \nabla g_t\big)\bigg)\on{vol}(g).
\end{aligned}}$$

The conditions of Lemma~\ref{nmb:3.4} are violated for general $g_0$.
However, we consider the special case that $M$ admits a metric $g_0$ such that the  Einstein equation 
$\on{Ricci}(g_0)= C g_0$ is satisfied. Let $g=r\g_0\in\NN$, then $\Scal(g)=\frac{Cn}{r}$.
For $h \in T_g(\NN)$  we calculate 
$$P_gh=\Phi(\Scal(g))h=\Phi(\frac{Cn}{r})h\in T_g(\NN)\;.$$
It remains to show that
$(D_{(g,m)}P)(g)$ is $\wt g^0_2$-normal to $\NN$ for all 
$m \in T_g\Met(M)$ such that $m$ is $\wt g^0_2$-normal to $\NN$. This follows
from  
\begin{align*}
 &\wt g^0_2\big((D_{(g,m)}P)g,g\big)= 
\\&
=\wt g^0_2\Big(\Phi'(\Scal(g))\big(\Delta(\Tr(g\i.m))+\na^*(\na^*(m))-g^0_2(\Ric,m)\big)g,g\Big)
\\&
=\Phi'\big(\frac{Cn}{r}\big) \int_M \big(\Delta(\Tr(g\i.m))+\na^*(\na^*(m))-
\\&\qquad\qquad\qquad\qquad\qquad\qquad\qquad
-g^0_2(\Ric,m)\big)g^0_2(g,g)\vol(g)
\\&
=\Phi'\big(\frac{Cn}{r}\big) n \int_M \big(\Delta(\Tr(g\i.m))+\na^*(\na^*(m))-g^0_2(Cg,m)\big)\vol(g)
=0\,;
\end{align*}
The first two terms vanish because they are divergences, and the last term vanishes by assumption 
on $m$.
Thus $\NN$ is a totally geodesic submanifold if $g_0$ satisfies the Einstein equation.
For the corresponding functions $\Psi$ and $f$ we obtain:
$$ \Psi(r):= \Phi(\frac1{r} \Scal(g_0))= \Phi(\frac{Cn}{r}), \qquad f(r)=- \Phi'(\frac{Cn}{r})\frac{Cn}{r^2}\;.$$ 
Thus $g(t)=r(t)g_0$ is a geodesic iff $g_0$ is a solution to the Einstein equation and $r$ satisfies
\begin{align*}
r'' \Phi\big(\frac{Cn}{r}\big) =\frac{ r'^2}{r} \Big(\frac12  \Phi'\big(\frac{Cn}{r}\big)\frac{Cn}{r}  +\left(1-\frac{n}{4}\right) \Phi\big(\frac{Cn}{r}\big)\Big) \enspace .
\end{align*}

In the case that $M$ does not admit a metric solving the Einstein equation  we cannot use Lemma~\ref{nmb:3.5} to check for geodesic incompleteness, but we can still compute the length of shrinking a metric to zero. 
Let $g(r)=rg_0$, with $\Scal(g_0)$ not necessary constant:
\begin{align*}
\on{Len}_0^1(g) &= 
\int_0^1 r^{\frac{n}{4}-1}\sqrt{n}\Big(\int_M \Ph\big(\tfrac{\Scal(g_0)}{r}\big) \vol(g_0)\Big)^{1/2}dr
\end{align*}
Now let us assume that $\Ph(u) \le C (1+|u|^{2k})$ for constants $C$ and $k$. 
\begin{align*}
&\on{Len}_0^1(g) \le
\int_0^1 r^{\frac{n}{4}-1}\sqrt{n}\Big(C\int_M \big(1+\tfrac{|\Scal(g_0)|^{2k}}{r^{2k}}\big) \vol(g_0)\Big)^{1/2}dr
\\&
=\int_0^1 r^{\frac{n}{4}-1}\sqrt{nC}\Big(\Vol(g_0)+\frac1{r^{2k}}\int_M |\Scal(g_0)|^{2k}\vol(g_0)\Big)^{1/2}dr
\\&
=\int_0^1 r^{\frac{n}{4}-1}\sqrt{nC\Vol(g_0)}\Big(1+\frac1{2r^{2k}\Vol(g_0)}\int_M |\Scal(g_0)|^{2k}\vol(g_0)\Big)dr.
\end{align*}
This is finite if and only if $\frac{n}{4}-1-2k>-1$, i.e., $n>8k$.
Thus $(\Met(M),G^{\Ph(\Scal)})$ is geodesically incomplete if $M$ admits a metric solving the 
Einstein equation and $\Ph(u)\le C(1+|u|^{2k})$ for 
$k<\dim(M)/8$.

\subsection{Sobolev metrics using the Laplacian}\label{nmb:4.4} 

We first consider the Sobolev metric of the form
$$G^{P}_g(h,k)=\int_M g^0_2\big((1+\De)^p h,k\big) \vol(g)$$
where $\De^g$ is the geometric Bochner-Laplacian described in \ref{nmb:2.8}. 
The adjoint of the derivative of $P$ satisfies
\begin{align*}
&\int_M g^0_2\big(m,(D_{(g,.)}Ph)^*(k)\big) \vol(g)=\int_M g^0_2\big((D_{(g,m)}P)h,k\big) \vol(g)\\&\qquad
=\sum_{i=1}^{p}\int_M g^0_2\big((1+\De)^{i-1}(D_{(g,m)}\De)(1+\De)^{p-i}h,k\big) \vol(g)\\&\qquad
=\sum_{i=1}^{p}\int_M g^0_2\big((D_{(g,m)}\De)(1+\De)^{p-i}h,(1+\De)^{i-1}k\big) \vol(g)\\&\qquad
=\sum_{i=1}^{p}\int_M g^0_2\Big(m,\big((D_{(g,.)}\De)(1+\De)^{p-i}h\big)^*(1+\De)^{i-1}k\Big) \vol(g) \enspace.
\end{align*}
Thus it remains to calculate the adjoint of the derivative of $\Delta$.
\begin{lem*}
The differential of the Laplacian acting on $\binom 02$-tensors  admits an adjoint with 
respect to the metric $\wt g_2^0$, which is given by:
\begin{multline*}
\wt g_2^0\big(D_{(g,m)} \De h ,k \big)=:\wt g_2^0\big(m ,(D_{(g,.)} \De h)^*(k) \big)=
\\
= \wt g_2^0\Big(\! m, g^{i_1j_1}g^{i_2j_2}(\na^2h)_{.. 
i_1i_2}k_{j_1j_2}-(N^0_3(.)\na h)^*(g\otimes k)+(N^0_2(.)h)^*(\na k)\! \Big)
\end{multline*}
Here $(N^0_q(.)h)^*$ denotes the adjoint of the differential of the covariant derivative:
$$\wt g^0_{q+1}(N^0_q(m)h, k)=: \wt g^0_2 \big(m, (N^0_q(.)h)^*(k)\big)
= \wt g^0_2 \big(m, \nabla^* (\si(N_q^0)(.)h)^*k\big),$$
where $h\in\Ga(T^0_{q}M)$, $k\in \Ga(T^0_{q+1}M)$ and where $\si(N_q^0)$ 
denotes the total symbol of $N_q^0$.
It is tensorial and of the form
\begin{align*}
&\si(N^0_q)(\wt m)(h)(X_0,\ldots, X_q)= 
\\&
= \frac{-1}2\sum_{j=1}^q h\Big(X_1,\ldots,X_{j-1},\sum_{i=0}^2(-1)^i 
  \big(\ta^i(\wt m)(X_0,X_j,\cdot)\big)^\sharp,X_{j+1},\ldots ,X_q\Big),
\end{align*}
where $\wt m \in \Ga(T^*M \otimes S^2T^*M)$, $h \in \Ga(T^0_qM)$, $X_0,\ldots,X_q \in \X(M)$, 
and where $\ta^i$ is the $i$-th power of the cyclic permutation $\ta(\al\otimes\be\otimes\ga)=\ga\otimes\al\otimes\be$. 
\end{lem*}
\begin{proof}
To prove the formula for $(N^0_q(.)h)^*$ it suffices to show that $$N^p_q(m)(h)=\si(N^p_q)(\na m)(h).$$
This follows from \eqref{eq:n10} and \eqref{eq:npq} in \ref{nmb:2.5}. The formula for $D_{(g,.)}\De$ follows 
from \ref{nmb:2.8}.
\end{proof}
The above discussion and the formula for the geodesic equation from Section~\ref{nmb:3.1}
yield the geodesic equation for Sobolev type metrics:
$$\boxed{\begin{aligned}
&((1+\De)^pg_{t})_t=\frac12 g^{i_1j_1}g^{i_2j_2}(\na^2(1+\De)^{p-i}g_t)_{\:\! .\, .\,  i_1i_2}(1+\De)^{i-1}(g_t)_{j_1j_2}
\\&\qquad
-\frac12(N^0_3(.)\na (1+\De)^{p-i}g_t)^*(g\otimes (1+\De)^{i-1}g_t)
\\&\qquad
+\frac12 (N^0_2(.)(1+\De)^{p-i}g_t)^*(\na (1+\De)^{i-1}g_t)
\\&\qquad
+\frac14.g.\Tr(g\i.(1+\De)^pg_t.g\i.g_t)
+ \frac12 g_t.g\i.(1+\De)^pg_t
\\&\qquad
+\frac12 (1+\De)^pg_t.g\i.g_t-\frac12\Tr(g\i.g_t).(1+\De)^pg_t.
\end{aligned}}$$
The conditions of Theorem~\ref{nmb:3.2} are valid, so the geodesic equation is well-posed. 
The conserved quantity is
$$\boxed{\begin{aligned}
 &\nabla^* \big((1+\De)^p g_t)\big)\on{vol}(g).
\end{aligned}}$$

Finally we want to study again the geodesics of pure scaling using Lemma~\ref{nmb:3.4}.  
Let $g_0 \in \Met(M)$ and $g =r g_0 \in \NN$. 
Since $\nabla g=0$ and consequently $\Delta g =0$, one has
$P_g g= (1+\Delta)^pg = g$. 
It remains to show that
$(D_{(g,m)}P)(g)$ is $\wt g^0_2$-normal to $\NN$ for all 
$m \in T_g\Met(M)$ such that $m$ is $\wt g^0_2$-normal to $\NN$. This follows
from
\begin{align*}
&\wt g^0_2\Big((D_{(g,m)}P)g,g\Big)=\sum_{i=1}^{p}\wt g^0_2\Big((1+\De)^{i-1}(D_{(g,m)}\De)(1+\De)^{p-i}g,g\Big)
\\&
=\sum_{i=1}^{p}\wt g^0_2\Big((D_{(g,m)}\De)(1+\De)^{p-i}g,(1+\De)^{i-1}g\Big)=
p\; \wt g^0_2\Big((D_{(g,m)}\De)g,g\Big)
\\&
=p\; \wt g^0_2\Big(\Tr(g\i m g\i \na^2 g) - \Tr^g\big(N^0_{3}(m) \na g\big) - \Tr^g\big(\na N^0_2(m) g\big),g\Big)
\\&
=0-0+p\; \wt g^0_2\Big( \na^*\big( N^0_2(m) g\big),g\Big)
=p\; \wt g^0_3\Big( \big( N^0_2(m) g\big), \na g\Big)=0.
\end{align*}
Thus the conditions of 
Lemma~\ref{nmb:3.4} are satisfied and $\NN$ is a totally geodesic submanifold for every $g_0\in\Met(M)$. 
Furthermore, since  $(D_{(rg_0,\cdot)}Pg_0)^*(g_0) = g_0$, the equation for geodesics of the form $r(t)g_0$ with respect to Sobolev metrics is equal 
to that with respect to the $L^2$ metric, c.f. Section~\ref{nmb:4.1}. In particular this proves that $\Met(M)$ is 
geodesically incomplete for each Sobolev metric.

\subsection{General Remarks}\label{nmb:4.5}
The $L^2$-metric is the only of the above discussed examples that it is relatively well-understood. 
An explicit analytic formula for geodesics has been derived, e.g. in \cite{Gil-MedranoMichor91}, 
and as a direct consequence it has been shown that
the space of Riemannian metrics is not complete with respect to this metric. 
Furthermore, the completion of this space has been described and analyzed in \cite{Clarke09a,Clarke09c}.

For the other metrics described in this section the situation is more complicated, 
since there is no hope to find  general analytic solutions to the corresponding geodesic equations.
But the equations as presented above are ready for numerical implementation.  This has been successfully done for the related
spaces of immersions and shapes, see \cite{Michor119,Michor118,MBMB2011}. 
Another issue is that we still do not know whether there exists a metric such that the  space of all Riemannian 
metrics is geodesically complete.

\section{The Ricci vector field}\label{nmb:5}

The space of metrics $\Met(M)$ is a convex open subset in the Fr\'echet space $\Ga(S^2T^*M)$. 
So it is contractible.  A necessary and sufficient condition for Ricci curvature
to be a gradient vector field with respect to the $G^P$-metric is that the following exterior derivative 
vanishes:
\begin{align*}
&\big( d G^P(\Ric, \cdot) \big) (h,k)  =
\\&\qquad
= h G^P(\Ric, k) - k G^P(\Ric, h) - G^P(\Ric,[h, k])=0.
\end{align*}
It suffices to look at constant vector fields $h,k$, in which case $[h,k]=0$. 
We have  
\begin{align*}
&h G^P(\Ric, k) - k G^P(\Ric, h)  
\\& 
= \int \Big( 
- \Tr\big(g\i h g\i  (P \Ric) g\i k\big) + \Tr\big(g\i k g\i  (P \Ric) g\i h\big) 
\\&\qquad\quad
+  \Tr\big(g\i D_{g,h} (P \Ric) g\i k\big) - \Tr\big(g\i D_{g,k} (P \Ric) g\i h\big) 
\\&\qquad\quad
-  \Tr\big(g\i  (P \Ric) g\i h g\i k\big) + \Tr\big(g\i  (P \Ric) g\i k g\i h\big) 
\\&\qquad\quad
+ \frac12  \Tr\big( g\i (P \Ric) g\i k \big) \Tr(g\i h) 
\\&\qquad\quad
- \frac12  \Tr\big( g\i (P \Ric) g\i h \big) \Tr(g\i k) \Big) \vol(g).
\end{align*}
Some terms in this formula cancel out because for symmetric $A,B,C$ one has 
$\Tr(ABC)=\Tr((ABC)^\top)=\Tr(C^\top B^\top A^\top)=\Tr(A^\top C^\top B^\top) = \Tr(ACB).$ Therefore
\begin{align*}
h G^P(&\Ric, k) - k G^P(\Ric, h)  \\ =
\int \Big( 
  & \Tr\big(g\i D_{g,h} (P \Ric) g\i k\big) - \Tr\big(g\i D_{g,k} (P \Ric) g\i h\big) \\
&+  \frac12 \Tr\big( g\i (P \Ric) g\i k \big) \Tr(g\i h) \\&- \frac12 \Tr\big( g\i (P \Ric) g\i h \big) \Tr(g\i k) 
\Big) \vol(g).
\end{align*}
We write
$D_{g,h} (P \Ric) = Q(h)$
for some differential operator $Q$ mapping symmetric two-tensors to themselves and $Q^*$ for the adjoint 
of $Q$ with respect to $\wt g^0_2$.
\begin{align*}
&h G^P(\Ric, k) - k G^P(\Ric, h)  =
\\& 
= \int \Big( g^0_2\big(Q(h), k\big) 
- g^0_2\big(Q(k), h\big) 
+ \frac12 g^0_2\big(P \Ric, k \big) \Tr(g\i h) 
\\&\qquad\qquad\qquad\qquad\qquad\qquad\qquad\;
- \frac12 g^0_2\big(P \Ric, h \big) \Tr(g\i k) \Big) \vol(g) 
\\&
= \int g^0_2\Big(Q(h)-Q^*(h)+\frac12 (P \Ric).\Tr(g\i h)
\\&\qquad\qquad\qquad\qquad\qquad\qquad\qquad\;
-\frac12 g.g^0_2\big(P \Ric, h\big),k\Big) \vol(g).
\end{align*}
We have proved:

\begin{lem*}
The Ricci vector field $\Ric$ is a gradient field for the $G^P$-metric if and only if the equation
\begin{equation}\begin{aligned}
2\big(Q(h)-Q^*(h)\big)+ (P \Ric).\Tr(g\i h)- g.g^0_2\big(P \Ric, h\big) = 0, \\
\text{with } Q(h)=Q_g(h)=D_{g,h} (P_g \Ric_g),
\end{aligned}\tag{1}\end{equation}
is satisfied for all $g\in\Met(M)$ and and all symmetric $\binom{0}{2}$-tensors $h$.
\end{lem*}

None of the specific metrics studied in Section~\ref{nmb:4} of this paper 
satisfies the Lemma in general dimension. 
Note that the Lemma is trivially satisfied in dimension $\on{dim}(M)=1$.
In dimension $2$ the equation $\Ric_g=\frac12 \Scal_g$ holds and the operator $P_gh=2\Scal_g\i h$
satisfies equation (1) on the open subset $\{g: \Scal_g \neq 0\}$.
Generally, equation (1) is satisfied if $P_g \Ric_g=g$, but this 
cannot hold on the space of all metrics if $\dim(M)>2$.

\end{document}